\documentclass[11pt,a4paper]{article}
\usepackage[latin1]{inputenc}
\usepackage{amsmath}
\usepackage{amsfonts}
\usepackage{amssymb}
\usepackage[dvips]{graphicx}
\usepackage{float}
\usepackage{stmaryrd}
\restylefloat{figure}

\newtheorem{satz}{Satz}[section]
\newtheorem{theorem}[satz]{Theorem}
\newtheorem{cor}[satz]{Corollary}
\newtheorem{defin}[satz]{Definition}

\newcommand{\abs}[1]{\left|{#1}\right|}
\newcommand{\norm}[1]{\left|\left|{#1}\right|\right|}
\newcommand{\rund}[1]{\left(#1\right)}
\newcommand{\spitz}[1]{\left\langle{#1}\right\rangle}
\newcommand{\schweif}[1]{\left\{#1\right\}}

\def\zz{\mathbb{Z}}
\def\cz{\mathbb{C}}
\def\nz{{\rm I\kern-.20em N}}
\def\rz{{\rm I\kern-.20em R}}

\def\B{\begin{cal} B \end{cal}}
\def\C{\begin{cal} C \end{cal}}
\def\D{\begin{cal} D \end{cal}}
\def\H{\begin{cal} H \end{cal}}
\def\O{\begin{cal} O \end{cal}}
\def\F{\begin{cal} F \end{cal}}

\def\g{\frak g}
\def\a{\frak a}
\def\n{\frak n}
\def\su{\frak su}

\def\vol{{\rm vol\ }}
\def\Im{{\rm Im\ }}
\def\tr{{\rm tr}}
\def\Leb{{\rm Leb}}
\def\cosh{{\rm cosh\ }}
\def\sinh{{\rm sinh\ }}
\def\Re{{\rm Re\ }}

\def\bz{{\bf z} }
\def\bv{{\bf v} }
\def\be{{\bf e} }
\def\bb{{\bf b} }
\def\bc{{\bf c} }
\def\bu{{\bf u} }
\def\bw{{\bf w} }
\def\bl{{\bf l} }

\def\b0{{\bf 0} }

\def\eps{\varepsilon }
\def\bzeta{{\bf \zeta} }
\def\bEta{{\bf \eta} }
\def\bthet{{\bf \vartheta} }

\begin{document}

\vskip 1.0 true cm

\begin{center}
{\Huge \bf A {\sc Satake} type theorem for

\vskip 0.4 true cm

Super Automorphic forms}
\end{center}

\vskip 1.0 true cm

\begin{center}
Roland Knevel ,

Philipps-Universität Marburg (Germany)
\end{center}

\vskip 2.0 true cm

\section*{Mathematical Subject Classification}

11F55 (Primary) , 32C11 (Secondary) .

\section*{Keywords}

Automorphic and cusp forms, complex bounded symmetric domains, super symmetry, semisimple {\sc Lie} groups, unbounded realization of a complex bounded symmetric domain.

\section*{Abstract}

Aim of this article is a {\sc Satake} type theorem for super automorphic forms on a complex bounded symmetric super domain $\B$ of rank $1$ with respect
to a lattice $\Gamma$ .  'Super' means: additional odd (anticommuting) coordinates on an ordinary complex bounded symmetric domain $B$ (the so-called body of $\B$ ) of rank $1$ . {\sc Satake}'s theorem says that for large weight $k$ all spaces

\[
sM_k(\Gamma) \cap L_k^s\rund{\Gamma \backslash \B} \, ,
\]

$s \in \, [ \, 1, \infty \, ] \, $ coincide, where $sM_k(\Gamma)$ denotes the space of super automorphic forms for $\Gamma$ with respect to the weight $k$ , and $L_k^s\rund{\Gamma \backslash \B}$ denotes the space of $s$-intergrable functions with respect 
to a certain measure on the quotient $\Gamma \backslash \B$ depending on $k$ . So all these spaces are equal to the space \\
$sS_k(\Gamma) := sM_k(\Gamma) \cap L_k^2\rund{\Gamma \backslash \B}$ of super cusp forms for $\Gamma$ to the weight $k$ . \\

As it is already well known for automorphic forms on ordinary complex bounded symmetric domains, we will give a proof of this theorem using an unbounded realization $\H$ of $\B$ and {\sc Fourier} decomposition at the cusps of the quotient
$\Gamma \backslash B$ mapped to $\infty$ via a partial {\sc Cayley} transformation.

\section*{Introduction}

Automorphic and cusp forms on an ordinary complex bounded symmetric domain $B$ are a classical field of research. Let us give a general definition:

\begin{defin}[automorphic and cusp forms in general] \label{autom forms general}
Suppose $B \subset \cz^n$ is a bounded symmetric domain and $G$ a semisimple {\sc Lie} group acting transitively and
holomorphically on $B$ . Let $j \in \C^\infty(G \times B, \cz)$ be a cocycle, this means $j$ is a smooth function on $G \times B$ , holomorphic in the second entry, such that

\[
j(g h, z) = j(g, h z) j(h, z)
\]

for all $g, h \in G$ and $z \in B$ . Let $k \in \zz$ and $\Gamma \sqsubset G$ be a
discrete subgroup.

\item[(i)] A holomorphic function $f \in \O(B)$ on $B$ is called an automorphic form of weight $k$ with respect to $\Gamma$ if and only if
$f = f|_\gamma$ for all~$\gamma~\in~\Gamma$~, where $f|_g \rund{\bz} := f\rund{g \bz} j\rund{g, \bz}^k$ for all $\bz \in B$ and
$g \in G$~, or equivalently the lift $\widetilde f \in \C^\infty(G)$ is left-$\Gamma$-invariant, where 
$\widetilde f (g) := f|_g \rund{\b0}$ for all $g \in G$ . The space of automorphic forms of weight $k$ with respect
to $\Gamma$ is denoted by $M_k(\Gamma)$~.

\item[(ii)] An automorphic form $f \in M_k(\Gamma)$ is called a cusp form of weight $k$ with respect to $\Gamma$ if and only if
$\widetilde f \in L^2\rund{\Gamma \backslash G}$ . The {\sc Hilbert} space of cusp forms of weight $k$ with respect to $\Gamma$ is denoted by $S_k(\Gamma)$ .
\end{defin}

In the simplest case, where $B \subset \cz$ is just the unit disc, $G = SU(1, 1)$ acting on $B$ via {\sc Möbius} transformations,

\[
j(g, z) = \frac{1}{c z + d} \, , \, g = \rund{\begin{array}{cc} a & b \\
c & d \end{array}} \in SU(1, 1) \, ,
\]

and $\Gamma \sqsubset G$ is a lattice, this means a discrete subgroup with finite covolume, one needs a more restrictive definition for automorphic and cusp forms. It is well known that after adding the cusps of $\Gamma \backslash B$ in $\partial B$ , 
which are always finitely many, the quotient $\Gamma \backslash B$ is compact. Having fixed a cusp $z_0 \in \partial B$ of $\Gamma \backslash B$ there exists a {\sc Cayley} transform $R$ mapping biholomorphically the unit disc $B$ onto the upper half 
plane $H \subset \cz$ and $z_0$ to $i \infty$ . Since $\Gamma$ is a lattice there exists an element $\gamma \in \Gamma$ such that

\[
R \gamma R^{- 1} = \rund{\begin{array}{cc} 1 & \lambda_0 \\
0 & 1 \end{array}} \, ,
\]

$\lambda_0 \setminus \{0\}$ , acting on $H$ as translation $w \mapsto w + \lambda_0$ . If a function $f \in \O(B)$ fulfills $f|_\gamma = f$ then $f|_{R^{- 1}} \in \O(H)$ fulfills

\[
f|_{R^{- 1}}(w) = \left.f|_{R^{- 1}} \right|_{R \gamma R^{- 1}} (w) = f|_{R^{- 1}}\rund{w + \lambda_0}  \, ,
\]

and so it has a {\sc Fourier} decomposition

\begin{equation}
f|_{R^{- 1}}(w) = \sum_{m \in \frac{1}{\lambda_0} \zz} c_m e^{2 \pi i m w}  \, . \label{Fourier 1dim}
\end{equation}

\begin{defin}[automorphic and cusp forms on the unit disc $B$ ] \label{autom forms onedim}

\item[(i)] A holomorphic function $f \in \O(B)$ is called an automorphic form of weight $k$ for $\Gamma$ if and only if $f|_\gamma = f$ for all $\gamma \in \Gamma$ and for each cusp $z_0 \in \partial B$ of $\Gamma \backslash B$ it has a
positive {\sc Fourier} decomposition, this means precisely $c_m = 0$ in (\ref{Fourier 1dim}) for all $m < 0$ , or equivalently $f|_{R^{- 1}} (w)$ is bounded for $\Im w \leadsto \infty$ .

\item[(ii)] An automorphic form $f \in M_k(\Gamma)$ is called a cusp form if and only if it has a {\bf strictly} positive {\sc Fourier} decomposition for each cusp $z_0 \in \partial B$ of $\Gamma \backslash B$~, which means $c_m = 0$ in (\ref{Fourier 1dim}) for
all $m \leq 0$ , or equivalently $f|_{R^{- 1}} (w) \leadsto 0$ for $\Im w \leadsto \infty$ .
\end{defin}

However, in contrast to the one dimensional case, for higher dimension \\
$n \geq 2$ , when $B \subset \cz^n$ is the unit ball, $G = SU(n, 1)$ acting on $B$ via {\sc Möbius} transformations,

\[
j(g, z) = \frac{1}{\bc \bz + d} \, , \, g = \rund{\begin{array}{c|c} A & \bb \\ \hline
\bc & d \end{array}} \begin{array}{l}
\} n \\
\leftarrow n + 1
\end{array} \in SU(n, 1) \, ,
\]

and $\Gamma \sqsubset G$ is a lattice, the situation is different: Then again one has partial {\sc Cayley} transforms $R$ mapping $B$ onto an unbounded realization $H$ of $B$ , which traditionally is a generalization of the right half plane instead of 
the upper half plane, but a holomorphic function $f \in \O(B)$ fulfilling $f|_\gamma = f$ for all $\gamma \in \Gamma$ automatically has a 'positive' {\sc Fourier} decomposition at each cusp, and therefore the general definition \ref{autom forms general} is
considered to be the right one. This is known as {\sc Köcher}'s principle, see for example in section 11.5 of \cite{Baily} . Futhermore {\sc Satake}'s theorem says that in this case for weight $k \geq 2 n$ all spaces

\[
M_k(\Gamma) \cap L_k^s\rund{\Gamma \backslash G}  \, ,
\]

$s \in \, [ \, 1, \infty \, ] \, $ , coincide, and therefore are equal to $S_k(\Gamma) = M_k(\Gamma) \cap L_k^2\rund{\Gamma \backslash G}$~, where

\[
L_k^s\rund{\Gamma \backslash G} := \schweif{f \in \cz^B \, \left| \, \widetilde f \in L^s\rund{\Gamma \backslash G}\right.} \, .
\]

The crucial argument is that for any function $f \in M_k(\Gamma)$ , $k \geq 2 n$ and $s \in \, [ \, 1, \infty \, ] \, $ the following are equivalent:

\begin{itemize}
\item[(i)] $f \in L_k^s\rund{\Gamma \backslash G}$
\item[(ii)] $f$ has a '{\bf strictly} positive' {\sc Fourier} decomposition at each cusp.
\end{itemize}

In \cite{Baily} one can find this theory in more generality. \\

Since in recent time super symmetry has become an important field of research for mathematics and physics, one is also interested in
super automorphic resp. super cusp forms on complex bounded symmetric super domains with even (commuting) and odd (anticommuting) coordinates, and this article generalizes {\sc Köcher}'s principle and {\sc Satake}'s theorem for super automorphic 
forms on the complex super unit ball $\B$ with the usual unit ball $B \in \cz^n$ , $n \geq 2$ , as body, see theorems \ref{Satake Fourier} (ii) and \ref{main} . \\

{\it Acknowledgement:} The present paper is part of my PhD thesis, so I would like to thank my doctoral advisor Professor H. {\sc Upmeier} for many
helpful comments and mentoring and all the other persons who accompanied me during the time I spent in Marburg.

\section{The general setting}

Let $n \in \nz$ , $n \geq 2$ , $r \in \nz$ and $\B := B^{n| r}$ be the unique complex $(n, r)$-dimensional super domain with the unit ball

\[
B := B^n := \schweif{\left. \bz \in \cz^n \, \right| \, \bz ^* \bz < 1} \subset \cz^n
\]

as body, holomorphic even (commuting) coordinate functions $z_1, \dots, z_n$ and holomorphic odd (anticommuting) coordinate functions $\zeta_1, \dots, \zeta_r$ . Let us denote the space of (smooth) super functions (with values in $\cz$ ) on $\B$ by 
$\D(\B)$ and the space of super holomorphic functions on $\B$ by $\O(\B) \sqsubset \D(\B)$~. Let $\wp(r) := \wp\rund{\{1, \dots, r\}}$ . Then one can decompose every $f \in \D(\B)$ uniquely as

\[
f = \sum_{I, J \in \wp(r)} f_{I J} \zeta^I \overline \zeta^J \, ,
\]

all $f_{I J} \in \C^\infty(B, \cz)$ , $I, J \in \wp(r)$ , where $\zeta^I := \zeta_{i_1} \cdots \zeta_{i_\rho}$ , \\
$I = \schweif{i_1, \dots, i_\rho} \in \wp(r)$ , $i_1 < \dots < i_\rho$ , and every $f \in \O(\B)$ uniquely as

\[
f = \sum_{I \in \wp(r)} f_I \zeta^I   \, ,
\]

where all $f_I \in \O(B)$ . So

\[
\D(\B) \simeq \C^\infty(B, \cz) \otimes \bigwedge\rund{\cz^r} \boxtimes \bigwedge\rund{\cz^r} = \C^\infty(B, \cz) \otimes \bigwedge\rund{\cz^{2 r}}
\]

and

\[
\O(\B) \simeq \O(B) \otimes \bigwedge\rund{\cz^r} \, .
\]

Define

\begin{eqnarray*}
&& G := sS\rund{U(n, 1) \times U(r)} \\
&& \phantom{12} := \schweif{\left.\rund{\begin{array}{c|c}
g' & 0 \\ \hline
0 & E
\end{array}} \in U(n, 1) \times U(r) \, \right| \, \det g' = \det E}  \, ,
\end{eqnarray*}

which is a real $\rund{(n + 1)^2 + r^2 - 1}$ -dimensional {\sc Lie} group. Then we have a holomorphic action of $G$ on $\B$ given by super fractional linear
({\sc Möbius}) transformations

\[
g \rund{\begin{array}{c} \bz \\ \hline
\bzeta \end{array}} := \rund{\begin{array}{c} \rund{A \bz + \bb} \rund{\bc \bz + d}^{- 1} \\ \hline
E \bzeta \rund{\bc \bz + d}^{- 1}
\end{array}} \, ,
\]

where we split

\[
g := \rund{\begin{array}{c|c}
\begin{array}{c|c}
A & \bb \\ \hline
\bc & d\end{array} & 0 \\ \hline
0 & E
\end{array}} \begin{array}{l}
\} n \\
\leftarrow n + 1 \\
\} r
\end{array} \, .
\]

The stabilizer subgroup of $\b0$ in $G$ is

\begin{eqnarray*}
&& K := sS\rund{\rund{U(n) \times U(1)} \times U(r)} \\
&& \, = \schweif{\left.\rund{\begin{array}{c|c}
\begin{array}{c|c}
A & 0 \\ \hline
0 & d
\end{array} & 0 \\ \hline
0 & E
\end{array}} \in U(n) \times U(1) \times U(r) \, \right| \, \det A \, d = \det E}  \, ,
\end{eqnarray*}

which is a maximal compact subgroup of $G$ . On $G \times B$ we define the cocycle $j \in \C^\infty(G \times B, \cz)$ as $j(g, \bz) := \rund{\bc \bz + d}^{- 1}$ for all $g \in G$ and $\bz \in \B$ . It is holomorphic in the second entry. Let $k \in \zz$ be fixed. 
Then we have a right-representation of $G$ on $\D(\B)$ given by

\[
|_g : \D(\B) \rightarrow \D(\B) \, , \, f|_g\rund{\begin{array}{c}
\bz \\ \hline
\bzeta
\end{array}} := f\rund{g \rund{\begin{array}{c}
\bz \\ \hline
\bzeta
\end{array}} } j(g, \bz)^k
\]

for all $g \in G$ , which is holomorphic, more precisely if $f \in \O(\B)$ and $g \in G$ then $f|_g \in \O(\B)$ . Finally let $\Gamma$ be a discrete subgroup of~$G$~.

\begin{defin}[super automorphic forms]
Let $f \in \O(\B)$ . Then $f$ is called a super automorphic form for $\Gamma$ of weight $k$ if and only if $f|_\gamma = f$ for all $\gamma \in \Gamma$~. We denote the space of super automorphic forms for $\Gamma$ of weight $k$ by $sM_k(\Gamma)$ .
\end{defin}

Let $\cz^{0| r}$ be the purely odd complex super domain with one point $\{0\}$ as body and odd coordinate functions $\eta_1, \dots, \eta_r$ . Then \\
$\D\rund{\cz^{0| r}} \simeq \bigwedge\rund{\cz^r} \boxtimes \bigwedge\rund{\cz^r} \simeq \bigwedge\rund{\cz^{2 r}}$ . Let us define a lift:

\begin{eqnarray*}
\widetilde{\phantom{1}} : \D(\B) &\rightarrow& \C^\infty(G, \cz) \otimes \D\rund{\cz^{0| r}} \simeq \C^\infty(G, \cz) \otimes \bigwedge\rund{\cz^r} \boxtimes \bigwedge\rund{\cz^r} \, , \\
f &\mapsto& \widetilde f  \, ,
\end{eqnarray*}

where

\[
\widetilde f(g) := f|_g\rund{\begin{array}{c}
\b0 \\ \hline
\bEta
\end{array}} = f\rund{g \rund{\begin{array}{c}
\b0 \\ \hline
\bEta
\end{array}} } j\rund{g, \b0}^k
\]

for all $f \in \D(\B)$ and $g \in G$ . Let $f \in \O(\B)$ . Then clearly \\
$\widetilde f \in \C^\infty(G, \cz) \otimes \O\rund{\cz^{0| r}}$ and $f \in sM_k(\Gamma) \Leftrightarrow \widetilde f \in \C^\infty\rund{\Gamma \backslash G, \cz} \otimes \O\rund{\cz^{0| r}}$ since for all $g \in G$

\[
\begin{array}{ccc}
\C^\infty(G) \otimes \D\rund{\cz^{0| r}} & \mathop{\longrightarrow}\limits^{l_g} & \C^\infty(G) \otimes \D\rund{\cz^{0| r}} \\
\uparrow_{\, \widetilde{\phantom{1}}} & \circlearrowleft & \uparrow_{\, \widetilde{\phantom{1}}} \\
\D(\B) & \mathop{\longrightarrow}\limits_{\phantom{1} |_g} & \D(\B)
\end{array}  \, ,
\]

where $l_g: \C^\infty(G) \otimes \D\rund{\cz^{0| r}} \rightarrow \C^\infty(G) \otimes \D\rund{\cz^{0| r}} \, , \, l_g(f)(h) := f(g h)$ simply denotes the left translation with the element $g \in G$ . \\

Let $\spitz{\phantom{1}, \phantom{1}}$ be the canonical scalar product on $\D\rund{\cz^{0| r}} \simeq \bigwedge\rund{\cz^{2 r}}$ (semi-linear in the second entry) . Then for all 
$a \in \D\rund{\cz^{0| r}}$ we write $\abs{a} := \sqrt{\spitz{a, a}}$ , and $\spitz{\phantom{1}, \phantom{1}}$ induces a 'scalar product'

\[
(f, h)_{\Gamma} := \int_{\Gamma \backslash G} \spitz{\widetilde h, \widetilde f}
\]

for all $f, g \in \D(\B)$ such that $\spitz{\widetilde h, \widetilde f} \in L^1(\Gamma \backslash G)$ and for all $s \in \, ] \, 0, \infty \, ]$ a 'norm'

\[
\norm{f}_{s, \Gamma} := \norm{\phantom{\frac{1}{1}} \abs{\widetilde f} \phantom{\frac{1}{1}}}_{s, \Gamma \backslash G}
\]

for all $f \in \D(\B)$ such that $\abs{\widetilde f} \in L^s\rund{\Gamma \backslash G}$ . Recall that the scalar product $(\phantom{1}, \phantom{1})_{\Gamma}$ and the norm $\norm{\phantom{1}}_{s, \Gamma}$ actually depend on the weight $k$ . Let us define

\[
L_k^s(\Gamma \backslash \B) := \schweif{f \in \D(\B) \, \left| \begin{array}{c} \\ \\ \end{array}
\widetilde f \in \C^\infty(\Gamma \backslash G, \cz) \otimes \D\rund{\cz^{0|r}} \, , \, \norm{f}_{s, \Gamma}^{(k)} < \infty\right.}
\]

for all $s \in \, ] \, 0, \infty \,  ]$ .

\begin{defin}[super cusp forms]
Let $f \in sM_k(\Gamma)$ . $f$ is called a super cusp form for $\Gamma$ of weight $k$ if and only if $f \in L_k^2(\Gamma \backslash \B)$ . The
$\cz$- vector space of all super cusp forms for $\Gamma$ of weight $k$ is denoted by $sS_k(\Gamma)$ . It is a {\sc Hilbert} space.
\end{defin}

Observe that $|_g$ respects the splitting

\[
\O(\B) = \bigoplus_{\rho = 0}^r \O^{(\rho)}(\B)
\]

for all $g \in G$ , where $\O^{(\rho)}(\B)$ is the space of all $f = \sum_{I \in \wp(r) \, , \, \abs{I} = \rho} f_I \zeta^I$ , all $f_I \in \O(\B)$ , $I \in \wp(r)$ , $\abs{I} = \rho$ , $\rho = 0, \dots, r$ , and $\, \widetilde{\phantom{1}} \,$ maps the space $\O^{(\rho)}(\B)$ into
$\C^\infty(G, \cz) \otimes \O^{(\rho)}\rund{\cz^{0| r}} \simeq \C^\infty(G, \cz) \bigwedge\nolimits^{(\rho)}\rund{\cz^r}$ . Therefore we have splittings

\[
sM_k(\Gamma) = \bigoplus_{\rho = 0}^r sM_k^{(\rho)}(\Gamma) \phantom{1} \text{ and } \phantom{1} sS_k(\Gamma) = \bigoplus_{\rho = 0}^r sS_k^{(\rho)}(\Gamma)  \, ,
\]

where $sM_k^{(\rho)}(\Gamma) := sM_k(\Gamma) \cap \O^{(\rho)}(\B)$ , $sS_k^{(\rho)}(\Gamma) := sS_k(\Gamma) \cap \O^{(\rho)}(\B)$ , $\rho = 0, \dots, r$ , and the last sum is orthogonal. \\

In the following we will use the {\sc Jordan} triple determinant $\Delta: \cz^n \times \cz^n \rightarrow \cz$ given by

\[
\Delta\rund{\bz, \bw} := 1 - \bw^* \bz
\]

for all $\bz, \bw \in \cz^n$ . Let us recall the basic properties:

\begin{itemize}
\item[(i)] $\abs{j\rund{g, \b0}} = \Delta\rund{g \b0, g \b0}^\frac{1}{2}$ for all $g \in G$ ,

\item[(ii)] $\Delta\rund{g \bz, g \bw} = \Delta\rund{\bz, \bw} j\rund{g, \bz} \overline{j\rund{g, \bw}}$ for all $g \in G$ and $\bz, \bw \in B$ , and

\item[(iii)] $\int_B \Delta\rund{\bz, \bz}^\lambda d V_{\Leb} < \infty$ if and only if $\lambda > - 1$ .
\end{itemize}

Since $\abs{\det\rund{\bz \mapsto g \bz}'} = \abs{j(g, \bz)}^{n + 1}$ and because of (i) we have the $G$-invariant volume element $\Delta(\bz, \bz)^{- (n + 1)} d V_{\Leb}$ on $B$ . \\

For all $I \in \wp(r)$ , $h \in \O(B)$ , $\bz \in B$ and $g = \rund{\begin{array}{c|c}
g' & 0 \\ \hline
0 & E
\end{array}} \in G$ we have

\[
\left.\rund{h \bzeta^I} \right|_g \rund{\bz} = h\rund{g' \bz} \rund{E \bEta}^I j\rund{g, \bz}^{k + \abs{I}} \, ,
\]

where $E \in U(r)$ . So for all $s \in \, ] \, 0, \infty \, ] \, $ , $f = \sum_{I \in \wp(r)} f_I \bzeta^I$ and \\
$h = \sum_{I \in \wp(r)} h_I \bzeta^I \in \O(\B)$ we obtain

\[
\norm{f}_{s, \Gamma} \equiv \norm{\, \sqrt{\sum_{I \in \wp(r)} \abs{f_I}^2 \Delta\rund{\bz, \bz}^{k + \abs{I}} } \, }_{s, \,  \Gamma \backslash B, \, \Delta\rund{\bz, \bz}^{- (n + 1)} d V_{\Leb}}  \, ,
\]

and

\[
(f, h)_\Gamma \equiv \sum_{I \in \wp(r)} \int_{\Gamma \backslash B} \overline{f_I} h_I \Delta\rund{\bz, \bz}^{k + \abs{I} - (n + 1)} d V_{\Leb}
\]

if $\spitz{\widetilde h , \widetilde f} \in L^1(\Gamma \backslash G)$ , where '$\equiv$' means equality up to a constant $\not= 0$ depending on $\Gamma$ , $k$ and $s$ .

\section{{\sc Satake}'s theorem in the super case}

We keep the notation of section 1 , in particular $n \in \nz$ , $n \geq 2$ . Here now the main theorem of the article, which is the analogon to {\sc Satake}'s theorem for super automorphic forms:

\begin{theorem} \label{main} Let $\rho \in \{0, \dots, r\}$ . Assume $\Gamma \sqsubset G$ is a lattice (discrete such that $\vol \rund{\Gamma \backslash G} < \infty$ , $\Gamma \backslash G$ not necessarily compact) . Then

\[
sS_k^{(\rho)}(\Gamma) = sM_k^{(\rho)}(\Gamma) \cap L_k^s \rund{\Gamma \backslash \B}
\]

for all $s \in \, [ \, 1, \infty \, ] \,$ and $k \geq 2 n - \rho$ .

\end{theorem}

If $\Gamma \backslash G$ is compact then the assertion is trivial. For the non-compact case we will give a proof in the end of this section using the so-called unbounded realization $\H$ of $\B$ , which we will develop in the following.

By the way, as for ordinary automorphic forms, theorem \ref{main} implies that $sS_k(\Gamma)$ is finite dimensional for $n \geq 2$ , $\Gamma \sqsubset G$ being a lattice and $k \geq 2 n$ via lemma 12 of \cite{Baily} section 10. 2 , which says the 
following: \\

\begin{quote}
Let $(X, \mu)$ be a locally compact measure space, where $\mu$ is a positive measure such that $\mu(X) < \infty$ . Let $\F$ be a closed subspace of $L^2(X, \mu)$ which is contained in
$L^\infty(X, \mu)$ . Then

\[
\dim \F < \infty \, .
\]

\end{quote}

{\bf From now on let $\Gamma \backslash G$ be not compact.} \\

Let $\g' = \su(n, 1)$ be the {\sc Lie} algebra of $G' := SU(n, 1)$ ,

\[
G' \hookrightarrow G \, , \, g' \mapsto \rund{\begin{array}{c|c} g' & 0 \\ \hline
0 & 1 \end{array}} \, ,
\]

and let $\a \sqsubset \g'$ be the standard {\sc Cartan} sub {\sc Lie} algebra of $\g'$ . Then $A := \exp_G \a$ is the common standard maximal split Abelian subgroup of $G'$ and $G$ , it is the image of the {\sc Lie} group embedding

\[
\rz \hookrightarrow G' \, , \, t \mapsto a_t := \rund{\begin{array}{c|c|c}
\cosh t & 0 & \sinh t \\ \hline
0 & 1 & 0 \\ \hline
\sinh t & 0 & \cosh t
\end{array}}
\begin{array}{l}
\leftarrow 1 \\
\rbrace n - 1 \\
\leftarrow n + 1
\end{array}  \, .
\]

Let $\n \sqsubset \g'$ be the standard maximal nilpotent sub {\sc Lie} algebra, which is at the same time the direct sum of all root spaces of $\g'$ of
positive roots with respect to $\a$ . Let $N := \exp \n$ . Then we have an {\sc Iwasawa} decomposition

\[
G = N A K  \, ,
\]

$N$ is $2$-step nilpotent, and so $N' := [N, N]$ is at the same time the center of~$N$~. \\

Now we transform the whole problem to the unbounded realization via the standard partial {\sc Cayley} transformation

\[
R := \rund{\begin{array}{c|c|c}
\frac{1}{\sqrt{2}} & 0 & \frac{1}{\sqrt{2}} \\ \hline
0 & 1 & 0 \\ \hline
- \frac{1}{\sqrt{2}} & 0 & \frac{1}{\sqrt{2}}
\end{array}}
\begin{array}{l}
\leftarrow 1 \\
\rbrace n - 1 \\
\leftarrow n + 1
\end{array} \in G'^\cz = SL(n + 1, \cz)
\]

mapping $B$ via {\sc Möbius} transformation biholomorphically onto the unbounded domain

\[
H := \schweif{\left.\bw = \rund{\begin{array}{c}
w_1 \\ \hline
\bw_2
\end{array}} \begin{array}{l}
\leftarrow 1 \\
\rbrace n - 1
\end{array} \in \cz^n \, \right| \, \Re w_1 > \frac{1}{2} \bw_2^* \bw_2}  \, ,
\]

which is a generalized right half plane, and $\be_1$ to $\infty$ . We see that

\[
R G' R^{- 1} \sqsubset G'^\cz = SL(n + 1, \cz) \hookrightarrow GL(n + 1, \cz) \times GL(r, \cz)
\]

acts holomorphically and transitively on $H$ via fractional linear transformations, and explicit calculations show that

\[
a'_t := R a_t R^{- 1} = \rund{\begin{array}{c|c|c}
e^t & 0 & 0 \\ \hline
0 & 1 & 0 \\ \hline
0 & 0 & e^{- t}
\end{array}}
\begin{array}{l}
\leftarrow 1 \\
\rbrace n - 1 \\
\leftarrow n + 1
\end{array}
\]

for all $t \in \rz$ , and $R N R^{- 1}$ is the image of

\[
\rz \times \cz^{n - 1} \rightarrow R G' R^{- 1} \, , \, \rund{\lambda, \bu} \mapsto n'_{\lambda, \bu} := \rund{\begin{array}{c|c|c}
1 & \bu^* & i \lambda + \frac{1}{2} \bu^* \bu \\ \hline
0 & 1 & \bu \\ \hline
0 & 0 & 1
\end{array}}  \, ,
\]

which is a smooth diffeomorphism onto its image, with the multiplication rule

\[
n'_{\lambda, \bu} n'_{\mu, \bv} = n'_{\lambda + \mu + \Im\rund{\bu^* \bv} , \bu + \bv}
\]

for all $\lambda, \mu \in \rz$ and $\bu, \bv \in \cz^{n - 1}$ and acting on $H$ as pseudo translations

\[
\bw \mapsto \rund{\begin{array}{c}
w_1 + \bu^* \bw_2 + i \lambda + \frac{1}{2} \bu^* \bu \\ \hline
\bw_2 + \bu
\end{array}}  \, .
\]

Define $j\rund{R, \bz} =  \frac{\sqrt{2}}{1 - z_1} \in \O(B)$ , \\
$j\rund{R^{- 1}, \bw} := j\rund{R, R^{- 1} \bw}^{- 1} = \frac{\sqrt{2}}{1 + w_1} \in \O(H)$ , and for all

\[
g = \rund{\begin{array}{c|c}
\begin{array}{c|c}
A & \bb \\ \hline
\bc & d
\end{array} & 0 \\ \hline
0 & E
\end{array}} \in R G R^{- 1}
\]

define

\[
j\rund{g, \bw} = j\rund{R, R^{- 1} g \bw} j\rund{R^{- 1} g R, R^{- 1} \bw} j\rund{R^{- 1}, \bw} = \frac{1}{\bc \bw + d} \, .
\]

Let $\H$ be the unique $(n, r)$-dimensional complex super domain with body $H$~, holomorphic even coordinate functions $w_1, \dots, w_n$ and holomorphic odd coordinate functions $\vartheta_1, \dots, \vartheta_r$ . $R$ commutes with all
$g \in Z\rund{G'}$ , where

\[
Z\rund{G'} = \schweif{\left.\rund{\begin{array}{c|c}
\eps \, 1 & 0 \\ \hline
0 & E
\end{array}} \begin{array}{l}
\} n + 1 \\
\} r
\end{array} \right| \, \eps \in U(1) , E \in U(r) , \eps^{n + 1} = \det E} \sqsubset K
\]

denotes the centralizer of $G'$ in $G$ , and we have a right-representation of the group $R G R^{- 1}$ on $\D(\H)$ given by

\[
|_g : \D(\H) \rightarrow \D(\H) \, , \, f|_g\rund{\begin{array}{c}
\bw \\ \hline
\bthet
\end{array}} := f\rund{g \rund{\begin{array}{c}
\bw \\ \hline
\bthet
\end{array}} } j\rund{g, \bw}^k
\]

for all $g \in R G R^{- 1}$ , which is again holomorphic. If we define

\[
|_R : \D(\H) \rightarrow \D(\B) \, , \, f|_R\rund{\begin{array}{c}
\bz \\ \hline
\bzeta
\end{array}} := f\rund{R \rund{\begin{array}{c}
\bz \\ \hline
\bzeta
\end{array}} } j\rund{R, \bz}^k
\]

and

\[
|_{R^{- 1}} : \D(\B) \rightarrow \D(\H) \, , \, f|_{R^{- 1}}\rund{\begin{array}{c}
\bw \\ \hline
\bthet
\end{array}} := f\rund{R^{- 1} \rund{\begin{array}{c}
\bw \\ \hline
\bthet
\end{array}} } j\rund{R^{- 1}, \bw}^k  \, ,
\]

then we see that again if $f \in \O(\H)$ then $f|_R \in \O(\B)$ , and if $f \in \O(\B)$ then $f|_{R^{- 1}} \in \O(\H)$ , and

\[
\begin{array}{ccc}
\phantom{1234} \D(\H) & \mathop{\longrightarrow}\limits^{\phantom{1} |_{R g R^{- 1}} } & \D(\H) \\
|_R \phantom{1} \downarrow & \circlearrowleft & \phantom{12} \downarrow \phantom{1} |_R \\
\phantom{1234} \D(\B) & \mathop{\longrightarrow}\limits_{\phantom{1} |_g} & \D(\B)
\end{array}   \, .
\]

Now define the {\sc Jordan} triple determinant $\Delta'$ on $H \times H$ , which is again holomorphic in the first and antiholomorphic in the second variable, as

\[
\Delta'\rund{\bz, \bw} := \Delta\rund{R^{- 1} \bz, R^{- 1} \bw} j\rund{R^{- 1}, \bz}^{- 1} \overline{j\rund{R^{- 1}, \bw}}^{- 1} = z_1 + \overline{w_1} - \bw_2^* \, \bz_2
\]

for all $\bz, \bw \in H$ . Clearly again $\abs{\det\rund{\bw \mapsto g \bw}'} = \abs{j\rund{g, \bw}}^{n + 1}$ and $\abs{j\rund{g, \be_1}} = \Delta'\rund{g \be_1, g \be_1}^{\frac{1}{2}}$ for all $g \in R G R^{- 1}$ , and so $\Delta'\rund{\bw, \bw}^{- (n + 1)} d V_\Leb$ is 
the $R G R^{- 1}$ -invariant volume element on $H$ . If $f = \sum_{I \in \wp(r)} f_I \bzeta^I \in \O(\B)$~, all $f_I \in \O(B)$ , $I \in \wp(r)$ , then

\[
f|_{R^{- 1}}\rund{\begin{array}{c}
\bw \\ \hline
\bthet
\end{array}} = \sum_{I \in \wp(r)} f_I \rund{R^{- 1} \bw} j\rund{R^{- 1}, \bw}^{k + \abs{I}} \bthet^I \in \O(\H)  \, ,
\]

and if $f = \sum_{I \in \wp(r)} f_I \bthet^I \in \O(\H)$ , all $f_I \in \O(H)$ , $I \in \wp(r)$ , and \\
$g = \rund{\begin{array}{c|c}
* & 0 \\ \hline
0 & E
\end{array}} \in R G R^{- 1}$ , $E \in U(r)$ , then

\[
f|_g\rund{\begin{array}{c}
\bw \\ \hline
\bthet
\end{array}} = \sum_{I \in \wp(r)} f_I \rund{g \bw} j\rund{g, \bw}^{k + \abs{I}} \rund{E \bthet}^I \in \O(\H)  \, .
\]

Let $\partial H = \schweif{\bw \in \cz^n \, \left| \, \Re w_1 = \frac{1}{2} \bw_2^* \bw\right.}$ be the boundary of $H$ in $\cz^n$ . Then $\Delta'$ and $\partial H$ are $R N R^{- 1}$ -invariant, and
$R N R^{- 1}$ acts transitively on $\partial H$ and on each

\[
\schweif{\bw \in H \, \left| \, \Delta'\rund{\bw, \bw} = e^{2 t}\right.} = R N a_t \b0  \, ,
\]

$t \in \rz$ .

\pagebreak

\begin{figure}[H]
\begin{center}
\includegraphics[width=0.9\textwidth]{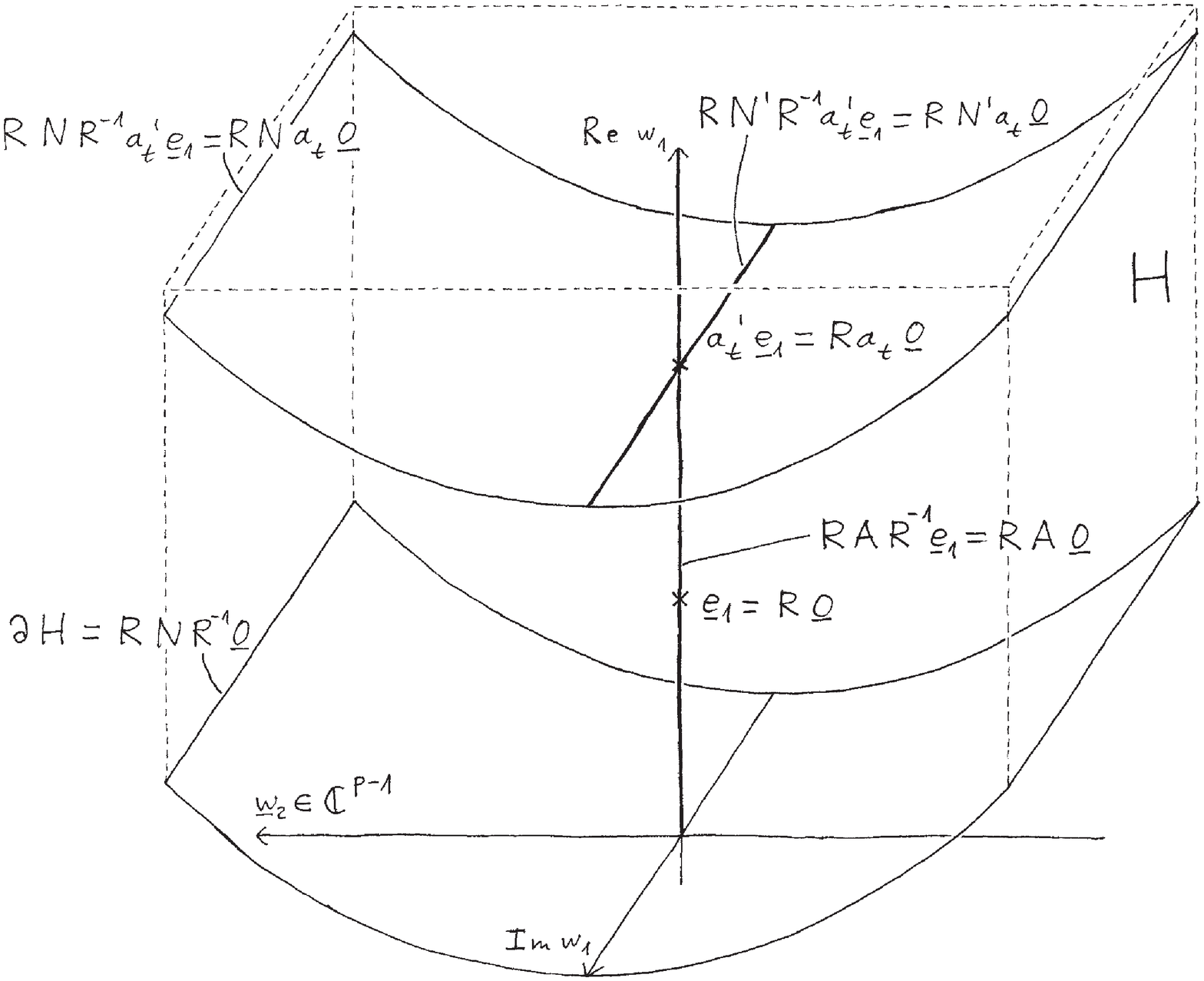}
\caption{the geometry of $H$ .}
\end{center}
\end{figure}

For all $t \in \rz$ define the rays $A_{< t} := \schweif{\left.a_\tau \, \right| \, \tau < t} \subset A$ and \\
$A_{> t} := \schweif{\left.a_\tau \, \right| \, \tau > t} \subset A$ .

\begin{theorem}[a 'fundamental domain' for $\Gamma \backslash G$ ] \label{fund} There exist $\eta \subset N$ open and relatively
compact , $t_0 \in \rz$ and $\Xi \subset G'$ finite such that if we define

\[
\Omega := \bigcup_{g \in \Xi} g \eta A_{> t_0} K
\]

then

\item[(i)] $g^{- 1} \Gamma g \cap N Z\rund{G'} \sqsubset N Z\rund{G'}$ and $g^{- 1} \Gamma g \cap N' Z\rund{G'} \sqsubset N' Z\rund{G'}$ are lattices, and

\[
N Z\rund{G'} = \rund{g^{- 1} \Gamma g \cap N Z\rund{G'}} \eta Z\rund{G'}
\]

for all $g \in \Xi$ ,
\item[(ii)] $G = \Gamma \Omega$ ,
\item[(iii)] the set $\schweif{\gamma \in \Gamma \, | \, \gamma \Omega \cap \Omega \not= \emptyset}$ is finite.
\end{theorem}

{\it Proof:} We use theorem 0.6 (i) - (iii) of \cite{GarlRagh} , which says the following:
\begin{quote}

Let $\Gamma' \subset G'$ be an admissible discrete subgroup of $G'$ . Then there exists $t'_0 > 0$ , an open, relatively compact subset $\eta_0 \subset N^+$~, a finite set $\Xi \subset G'$ , and
an open, relatively compact subset $\Omega'$ of $G'$ ( $\Xi$ being empty if $G' / \Gamma'$ is compact, and $\Omega'$ being empty if $G' / \Gamma'$ is non-compact) such that

\begin{itemize}
\item[(i)] For all $b \in \Xi$ , $\Gamma \cap b^{- 1} N^+ b$ is a lattice in $b^{- 1} N^+ b$ .
\item[(ii)] For all $t > t'_0$ and for all open, relatively compact subsets $\eta$ of $N^+$ such that $\eta \supset \eta_0$ , if

\[
\Omega'_{t, \eta} = \Omega' \cup \rund{\bigcup_{b \in \Xi} \sigma_{t, \eta} b}  \, ,
\]

then $\Omega'_{t, \eta} \Gamma' = G'$ , and
\item[(iii)] the set $\schweif{\gamma' \in \Gamma' \, , \, \Omega'_{t, \eta} \gamma' \cap \Omega'_{t, \eta} \not= \emptyset}$ is finite.
\end{itemize}

\end{quote}

Hereby $G'$ is a connected semisimple {\sc Lie} group of real rank $1$ , $N^+ = N$ is the standard nilpotent sub {\sc Lie} group of $G'$ and $\sigma_{t, \eta} := K' A_{< t} \eta$ for all $t > 0$ and
$\eta \subset N^+$ open and relatively compact, where $A$ denotes the standard maximal non-compact abelian and $K'$ the standard maximal compact sub {\sc Lie} group of $G'$ . Admissibility
is a geometric property of the quotient $\Gamma' \backslash G' / K'$ , roughly speaking $\Gamma'$ is called admissible if and only if $\Gamma' \backslash G' / K'$ has only finitely many
cusps. \\

Let us apply theorem 0.6 (i) - (iii) of \cite{GarlRagh} with $G' = SU(n, 1) \hookrightarrow G$ , \\
$K' := K \cap G' = S\rund{U(n) \times U(1)}$ and

\[
\Gamma' := \schweif{\left.\gamma' \in G' \, \right| \, \text{  there exists  } w \in Z\rund{G'} \text{  such that  } \gamma' w \in \Gamma} \sqsubset G'  \, ,
\]

which is of course again a lattice such that $\Gamma' \backslash G'$ is not compact and so it is admissible in the sense of \cite{GarlRagh} by theorem 0.7
of \cite{GarlRagh} . By lemma 3.18 of \cite{GarlRagh} $g^{- 1} \Gamma' g \cap N' \sqsubset N'$ is a lattice, and lemma 3.16 of \cite{GarlRagh} applied
with any $\rho \in \Gamma' \cap N' \setminus \{1\}$ tells us that $\left.\rund{g^{- 1} \Gamma' g \cap N} \right\backslash N$ is compact. So we see that there
exist $t_0 \in \rz$ , $\eta \subset N$ open and relatively compact and $\Xi \subset G'$ finite such that for all $g \in \Xi$

\[
\Gamma' \cap g N g^{- 1} \sqsubset g N g^{- 1}
\]

is a lattice, $\Gamma' \Omega' = G'$ if we define $\Omega' = \bigcup_{b \in \Xi} b \eta A_{< t_0} K'$ and

\[
\Delta := \schweif{\gamma' \in \Gamma' \, \left| \, \gamma' \Omega' \cap \Omega' \not= \emptyset\right.}
\]

is finite.

(i) and (ii) : now trivial by definition of $\Gamma' \sqsubset G'$ . $\Box$ \\

(iii) : Let $\gamma = \gamma' w \in \Gamma$ , $\gamma' \in \Gamma'$ , $w \in Z\rund{G'}$ , such that $\gamma \Omega \cap \Omega \not= \emptyset$ . Then

\[
\gamma' \Omega' Z\rund{G'} \cap \Omega' Z\rund{G'} \not= \emptyset  \, .
\]

Since $Z\rund{G'} \cap G' \sqsubset K'$ and $\Omega'$ is right-$K'$-invariant we have $\gamma' \Omega' \cap \Omega' \not= \emptyset$ as well and therefore $\gamma' \in \Delta$ . Conversely
$\gamma' Z\rund{G'}$ is compact and therefore $\Gamma \cap \gamma' Z\rund{G'}$ is finite for all $\gamma' \in \Gamma'$ . $\Box$ \\

From the 'fundamental domain $\Omega := \bigcup_{g \in \Xi} g \eta A_{> t_0} K$ one can really deduce the position of the cusps of $\Gamma \backslash B$ in $\partial B$ : they are up to the action of $\Gamma$ on $\partial B$ the limit points

\[
\lim_{t \rightarrow + \infty} g a_t \b0 = g \be_1 \, ,
\]

$g \in \Xi$ , where the limits are taken with respect to the Euclidian metric on $\cz^n$~. Their number is bounded above by $\abs{\Xi}$ and is therefore finite, as expected.

\begin{cor} \label{L^s} Let $t_0 \in \rz$ , $\eta \subset N$ and $\Xi \subset G$ be given by theorem
\ref{fund}~. Let $h \in \C\rund{\Gamma \backslash G, \cz}$ and $s \in \, ] \, 0, \infty \, ] \,$ . Then $h \in L^s\rund{\Gamma \backslash G}$ if and only if
$h\rund{g \bw} \in L^s\rund{\eta A_{> t_0} K}$ for all $g \in \Xi$ .
\end{cor}

{\it Proof:} If $s = \infty$ then it is evident since $G = \Gamma \Omega$ by theorem \ref{fund} (ii) . Now assume $s \in \, ] \, 0, \infty \, [ \,$ and $h \in L^s\rund{\Gamma \backslash G}$ .

\[
S := \abs{\schweif{\gamma \in \Gamma \, | \, \gamma \Omega \cap \Omega \not= \emptyset}} < \infty
\]

by theorem \ref{fund} (iii) . So for all $g \in \Xi$ we have

\[
\int_{\eta A_{> t_0} K} \abs{h\rund{g \bw}}^s = \int_{g \eta A_{> t_0} K} \abs{h}^s \leq \int_{\Omega} \abs{h}^s \leq S \int_{\Gamma \backslash G} \abs{h}^s < \infty  \, .
\]

Conversely assume $h\rund{g \bw} \in L^s\rund{\eta A_{> t_0} K}$ for all $g \in \Xi$ . Then since $G = \Gamma \Omega$ by theorem \ref{fund} (ii) we obtain

\[
\int_{\Gamma \backslash G} \abs{h}^s \leq \int_{\Omega} \abs{h}^s \leq \sum_{g \in \Xi} \int_{\eta A_{> t_0} K} \abs{h\rund{g \bw}}^s < \infty  \, . \, \Box
\] \\

Let $f \in sM_k(\Gamma)$ and $g \in \Xi$ . Then we may decompose

\[
\left.f|_g \right|_{R^{- 1}} = \sum_{I \in \wp(r)} q_I \bthet^I \in \O(\H) \, ,
\]

all $q_I \in \O(H)$ , $I \in \wp(r)$~, and by theorem \ref{fund} (i) we know that \\
$g^{- 1} \Gamma g \cap N' Z\rund{G'} \not\sqsubset Z\rund{G'}$ . So let $n \in g^{- 1} \Gamma g \cap N' Z\rund{G'} \setminus Z\rund{G'}$ ,

\[
R n R^{- 1} = n'_{\lambda_0, \b0} \rund{\begin{array}{c|c}
\eps 1 & 0 \\ \hline
0 & E
\end{array}}  \, ,
\]

$\lambda_0 \in \rz \setminus \{0\}$ , $\eps \in U(1)$ , $E \in U(r)$ , $\eps^{n + 1} = \det E$ . \\
$j\rund{R n R^{- 1}} := j\rund{R n R^{- 1}, \bw} = \eps^{- 1} \in U(1)$ is independent of $\bw \in H$ . So there exists $\chi \in \rz$ such that
$j\rund{R n R^{- 1}} = e^{2 \pi i \chi}$ . Without loss of generality we can assume that $E$ is diagonal, otherwise conjugate $n$ with an appropriate
element of $Z\rund{G'}$ . So there exists $D \in \rz^{r \times r}$ diagonal such that $E = \exp\rund{2 \pi i D}$ . If $D = \rund{\begin{array}{ccc}
d_1 & & 0 \\
 & \ddots & \\
0 & & d_r
\end{array}}$ and $I \in \wp(r)$ then we define $\tr_I \, D := \sum_{j \in I} d_j$ .

\begin{theorem}[{\sc Fourier} expansion of $\left.f|_g \right|_{R^{- 1}}$ ] \label{Satake Fourier}

\item[(i)] There exist unique $c_{I, m} \in \O\rund{\cz^{n - 1}}$ , $I \in \wp(r)$ , \\
$m \in \frac{1}{\lambda_0} \rund{\zz - \tr_I D - \rund{k + \abs{I}} \chi}$ , such that

\[
q_I\rund{\bw} = \sum_{m \in \frac{1}{\lambda_0} \rund{\zz - \tr_I D - \rund{k + \abs{I}} \chi}} c_{I, m}\rund{\bw_2} e^{2 \pi m w_1}
\]

for all $\bw \in H$ and $I \in \wp(r)$ , and so

\[
\left.f|_g \right|_{R^{- 1}}\rund{\bw} = \sum_{I \in \wp(r)} \phantom{1} 
\sum_{m \in \frac{1}{\lambda_0} \rund{\zz - \tr_I D - \rund{k + \abs{I}} \chi}} c_{I, m}\rund{\bw_2} e^{2 \pi m w_1} \bthet^I
\]

for all $\bw = \rund{\begin{array}{c}
w_1 \\ \hline
\bw_2
\end{array}} \begin{array}{l}
\leftarrow 1 \\
\rbrace n - 1
\end{array} \in H$ , where the convergence is absolute and compact.

\item[(ii)] $c_{I, m} = 0$ for all $I \in \wp(r)$ and $m > 0$  , and if \\
$\tr_I D + \rund{k + \abs{I}} \chi \equiv 0 \mod \zz$ in the group $(\rz, +)$ then $c_{I, 0}$ is a constant.

\begin{quote}
{\bf This is the super analogon for {\sc Köcher}'s principle, see section 11.5 of \cite{Baily} .} The condition $m > 0$ instead of $m < 0$ in definition \ref{autom forms onedim} comes from the fact that $\H$ is a generalized right half plane instead of the upper 
half plane.
\end{quote}

\item[(iii)] Let $I \in \wp(r)$ and $s \in \, [ \, 1, \infty \, ] \,$ . If $\tr_I D + \rund{k + \abs{I}} \chi \not\equiv 0 \mod \zz$ then

\[
q_I \Delta'\rund{\bw, \bw}^{\frac{k + \abs{I}}{2}} \in L^s\rund{R \eta A_{> t_0} \b0}
\]

with respect to the $R G R^{- 1}$ -invariant measure $\Delta'\rund{\bw, \bw}^{- (n + 1)} d V_\Leb$ on $H$~. If $\tr_I D + \rund{k + \abs{I}} \chi \equiv 0 \mod \zz$ and $k \geq 2 n - \abs{I}$ then

\[
q_I \Delta'\rund{\bw, \bw}^{\frac{k + \abs{I}}{2}} \in L^s\rund{R \eta A_{> t_0} \b0}
\]

with respect to the $R G R^{- 1}$ -invariant measure on $H$ if and only if $c_{I, 0} = 0$ .
\end{theorem}

{\it Proof:} (i) $f|_g$ is $g^{- 1} \Gamma g$ invariant, so we see that for all $\bw \in H$

\begin{eqnarray*}
\sum_{I \in \wp(r)} q_I\rund{\bw} \bthet^I &=& \left.f|_g \right|_{R^{- 1}} \rund{\bw} \\
&=& \left.\left.f|_g \right|_n \right|_{R^{- 1}} \rund{\bw} \\
&=& \sum_{I \in \wp(r)} q_I\rund{\bw + i \lambda_0 \be_1} \rund{E \bthet j\rund{R n R^{- 1}} }^I j\rund{R n R^{- 1}}^k \\
&=& \sum_{I \in \wp(r)} q_I\rund{\bw + i \lambda_0 \be_1} e^{2 \pi i \rund{\tr_I D + \rund{k + \abs{I}} \chi}} \bthet^I  \, .
\end{eqnarray*}

Therefore for all $\bw \in H$ and $I \in \wp(r)$

\[
q_I\rund{\bw} = q_I\rund{\bw + i \lambda_0 \be_1} e^{2 \pi i \rund{\tr_I D + \rund{k + \abs{I}} \chi}}  \, .
\]

Let $I \in \wp(r)$ . Then $h \in \O(H)$ given by

\[
h\rund{\bw} := q_I\rund{\bw} e^{-2 \pi i \frac{1}{\lambda_0} \rund{\tr_I D + \rund{k + \abs{I}} \chi} w_1}
\]

for all $\bw \in H$ is $i \lambda_0 \be_1$ periodic, and therefore there exists $\widehat h$ holomorphic on

\[
\widehat H := \schweif{\left.\bz = \rund{\begin{array}{c}
z_1 \\ \hline
\bz_2
\end{array}} \begin{array}{l}
\leftarrow 1 \\
\rbrace n - 1
\end{array} \, \right| \abs{z_1} > e^{\frac{\pi}{\lambda_0} \bz_2^* \bz_2}}
\]

such that for all $\bw \in H$

\[
h\rund{\bw} = \widehat h \rund{\begin{array}{c}
e^{\frac{2 \pi}{\lambda}} w_1 \\ \hline
\bw_2
\end{array}}   \, .
\]

{\sc Laurent} expansion now tells us that there exist $a_{m', \bl} \in \cz$ , $m' \in \zz$ , $\bl \in \nz^{n - 1}$ , such that

\[
\widehat h \rund{\bz} = \sum_{m' \in \zz} \, \sum_{\bl \in \nz^{n - 1}} a_{m', \bl} z_1^{m'} \bz_2^{\bl}
\]

for all $\bz = \rund{\begin{array}{c}
z_1 \\ \hline
\bz_2
\end{array}} \begin{array}{l}
\leftarrow 1 \\
\rbrace n - 1
\end{array} \in \widehat H$ , where the convergence is absolute and compact. Now let us define $d_{m'} \in \O\rund{\cz^{n - 1}}$ as

\[
d_{m'} \rund{\bz} := \sum_{\bl \in \nz^{n - 1}} a_{m', \bl} \bz_2^{\bl}  \, ,
\]

$m' \in \zz$ . Then for all $\bw \in H$

\[
q_I\rund{\bw} e^{- \frac{2 \pi i}{\lambda_0} \rund{\tr_I D + \rund{k + \abs{I}} \chi} w_1} = h\rund{\bw}
= \sum_{m' \in \zz} d_{m'}\rund{\bw_2} e^{\frac{2 \pi}{\lambda_0} m' w_1}  \, .
\]

So taking $c_m := d_{\lambda_0 m + \tr_I D + \rund{k + \abs{I}} \chi}$ , $m \in \frac{1}{\lambda_0} \rund{\zz - \tr_I D - \rund{k + \abs{I}} \chi}$ , gives the desired result. Uniqueness follows from
standard {\sc Fourier} theory. $\Box$ \\

(ii) {\it Step I } {\bf Show that all $q_I$ , $I \in \wp(r)$ , are bounded on $R N \b0 = \schweif{\left. \bw \in H \, \right| \, \Delta'\rund{\bw, \bw} = 2}$ .} \\

Obviously all $q_I$ , $I \in \wp(r)$ , are bounded on $R \eta \b0$ since $R \eta \b0$ lies relatively compact in $H$ . Let $C \geq 0$ such that $\abs{q_I} \leq C$ on $R \eta \b0$ for all $I \in \wp(r)$ .
By theorem \ref{fund}

\[
R N \b0 = R \rund{g^{- 1} \Gamma g \cap N Z\rund{G'}} \eta \b0 \,  .
\]

So let $R n' R^{- 1} = n'_{\lambda', \bu} \rund{\begin{array}{c|c}
\eps' 1 & 0 \\ \hline
0 & E'
\end{array}} \in g^{- 1} \Gamma g \cap N Z\rund{G'}$ , $\lambda' \in \rz$ , $\bu \in \cz^{n - 1}$~, $\eps' \in U(1)$ and $E' \in (r)$ . Then again

\[
j\rund{R n' R^{- 1}} := j\rund{R n' R^{- 1}, \bw} = \eps'^{- 1} \in U(1)
\]

is independent of $\bw \in H$ . Now if we use that $f \in sM_k(\Gamma)$ we get

\begin{eqnarray*}
\sum_{I \in \wp} q_I \bthet^I &=& \left.f|_g \right|_{R^{- 1}} \\
&=& \left.\left.f|_g \right|_{n'} \right|_{R^{- 1}} \\
&=& \sum_{I \in \wp(r)} q_I\rund{R n' R^{- 1} \bw} \rund{E' \bthet}^I \eps'^{k + \abs{I}}   \, .
\end{eqnarray*}

$\bigwedge\rund{\cz^r} \rightarrow \bigwedge\rund{\cz^r} \, , \, \bthet^I \mapsto \rund{E' \bthet}^I \eps'^{k + \abs{I}}$ is unitary, therefore

\[
\abs{q_I} \leq 2^r \abs{q_I\rund{R n' R^{- 1} \bw} }  \, .
\]

We see that $\abs{q_I} \leq 2^r C$ on the whole $R N \b0$ . \\

{\it Step II } {\bf Show that

\[
\abs{c_{I, m}\rund{\bw_2} e^{2 \pi m w_1}} \leq \norm{q_I}_{\infty, R N \b0}
\]

on $R N \b0$ for all $I \in \wp(r)$ and $m \in \frac{1}{\lambda_0} \rund{\zz - \tr_I D - \rund{k + \abs{I}} \chi}$ .} \\

Let $I \in \wp(r)$ and $m \in \frac{1}{\lambda_0} \rund{\zz - \tr_I D - \rund{k + \abs{I}} \chi}$ . By classical {\sc Fourier} analysis

\[
c_{I, m}\rund{\bw_2} e^{2 \pi m w_1} = \frac{1}{\lambda_0} \int_0^{\lambda_0} q_I\rund{\bw + i \lambda \be_1} e^{- 2 \pi i m \lambda} d \lambda
\]

for all $\bw \in H$ , and since $\bw + i \lambda \be_1 = n'_{\lambda, \b0} \bw \in R N R^{- 1} \bw$ the claim follows. \\

{\it Step III } {\bf Conclusion.} \\

Let $I \in \wp(r)$ and $m \in \frac{1}{\lambda_0} \rund{\zz - \tr_I D - \rund{k + \abs{I}} \chi}$ . Let $\bu \in \cz^{n - 1}$ be arbitrary. Then

\[
\rund{\begin{array}{c}
1 + \frac{1}{2} \bu^* \bu \\ \hline
\bu
\end{array}} \in R N \b0   \, ,
\]

and so

\[
\abs{c_{I, m}\rund{\bu}} \leq \norm{q_I}_{\infty, R N \b0} e^{- \pi m \bu^* \bu}  \, .
\]

Now the assertion follows by {\sc Liouville}'s theorem, {\bf where $n \geq 2$ is of course essential}. $\Box$ \\

(iii) Let

\[
\eta' := \schweif{\rund{i y, \bu} \in i \rz \oplus \cz^{n - 1} \, \left| \, \rund{\begin{array}{c}
1 + \frac{1}{2} \bu^* \bu + i y \\ \hline
\bu
\end{array}} \in R \eta \b0\right.}
\]

be the projection of $R \eta \b0$ onto $i \rz \oplus \cz^{n - 1}$ in direction of $\Re w_1 \in \rz$ . Then

\[
\Psi: \rz_{> e^{2 t_0}} \times \eta' \rightarrow R \eta A_{> t_0} \b0 \, , \, \rund{x, i y, \bu} \mapsto \rund{\begin{array}{c}
x + \frac{1}{2} \bu^* \bu + i y \\ \hline
\bu
\end{array}}
\]

is a $\C^\infty$-diffeomorphism with determinant $1$ , and

\[
\Delta'\rund{\Psi\rund{x, i y, \bu}, \Psi\rund{x, i y, \bu}} = 2 x
\]

for all $\rund{x, i y, \bu} \in \rz_{> e^{2 t_0}} \times \eta'$ . So

\[
q_I \Delta'\rund{\bw, \bw}^{\frac{k + \abs{I}}{2}} \in L^s\rund{R \eta A_{> t_0} \b0}
\]

with respect to the measure $\Delta'\rund{\bw, \bw}^{- (n + 1)} d V_\Leb$ if and only if

\[
\rund{q_I \circ \Psi} x^{\frac{k + \abs{I}}{2}} \in L^s\rund{\rz_{> e^{2 t_0}} \times \eta'}
\]

with respect to the measure $x^{- (n + 1)} d V_\Leb$ . \\

Now assume either $\tr_I D + \rund{k + \abs{I}} \chi \not\equiv 0 \mod \zz$ or \\
$\tr_I D + \rund{k + \abs{I}} \chi \equiv 0 \mod \zz$ and $c_{I, 0} = 0$ . Then in both cases by (ii) we can write

\[
q_I\rund{\bw} = \sum_{m \in \frac{1}{\lambda_0} \rund{\zz - \tr_I D - \rund{k + \abs{I}} \chi} \cap \rz_{< 0}} c_{I, m}\rund{\bw_2} e^{2 \pi m w_1}
\]

for all $\bw \in H$ , where the sum converges absolutely and uniformly on compact subsets of $H$ . Let us define

\[
M_0 := \max \frac{1}{\lambda_0} \rund{\zz - \tr_I D - \rund{k + \abs{I}} \chi} \cap \rz_{< 0} < 0 \, .
\]

Then since $R \eta a_{t_0} \b0\subset H$ is relatively compact and the {\sc Fourier} expansion in (i) has compact convergence we can define

\begin{eqnarray*}
&& C'' := e^{- 2 \pi M_0 e^{2 t_0}} \sum_{m \in \frac{1}{\lambda_0} \rund{\zz - \tr_I D - \rund{k + \abs{I}} \chi} \cap \rz_{< 0}} \norm{c_{I, m}\rund{\bw_2} e^{2 \pi m w_1}}_{\infty, R \eta a_{t_0} \b0} \\
&& \phantom{12345678901234567890123456789012345678901234567890} < \infty  \, .
\end{eqnarray*}

So we see that

\[
\abs{q_I\rund{\bw}} \leq C'' e^{\pi M_0 \Delta'\rund{\bw, \bw}}
\]

for all $\bw \in R \eta A_{> t_0} \b0$ ,

\[
\abs{q_I \circ \Psi} \leq C'' e^{2 \pi M_0 x}  \, ,
\]

and so $x^{\frac{k + \abs{I}}{2}} \rund{q_I \circ \Psi} \in L^s\rund{\rz_{> e^{2 t_0}} \times \eta'}$ with respect to the measure $x^{- (n + 1)} d V_\Leb$ . \\

Conversely assume $\tr_I D + \rund{k + \abs{I}} \chi \equiv 0 \mod \zz$ , $k \geq 2 n - \abs{I}$ and $c_{I, 0} \not= 0$~. Then as before we have the estimate

\[
\abs{\sum_{m \in \frac{1}{\lambda_0} \zz_{< 0}} c_{I, m}\rund{\bw_2} e^{2 \pi m w_1}} \leq C'' e^{- \pi \Delta'\rund{\bw, \bw}}
\]

for all $\bw \in R \eta A_{> t_0} \b0$ if we define

\[
C'' := e^{2 \pi e^{2 t_0}} \sum_{m \in \frac{1}{\lambda_0} \zz_{< 0}} \norm{c_{I, m}\rund{\bw_2} e^{2 \pi m w_1}}_{\infty, R \eta a_{t_0} \b0} < \infty  \, .
\]

Therefore there exists $S \geq 0$ such that

\[
\abs{\sum_{m \in \frac{1}{\lambda_0} \zz_{< 0}} c_{I, m}\rund{\bw_2} e^{2 \pi m w_1}} \leq \frac{1}{2} \abs{c_{I, 0}}  \, ,
\]

and so $\abs{q_I\rund{\bw}} \geq \frac{1}{2} \abs{c_{I, 0}}$ for all $\bw \in R \eta A_{> t_0} \b0$ having $\Delta'\rund{\bw, \bw} \geq S$~. So
$\abs{\rund{q_I \circ \Phi}\rund{x, i y, \bu}} \geq \frac{1}{2} \abs{c_{I 0}}$ for all $\rund{x, i y, \bu} \in \rz_{\geq S} \times \eta'$~, and so definitely
$x^{\frac{k + \abs{I}}{2}} \rund{q_I \circ \Phi} \notin L^s\rund{\rz_{> e^{2 t_0}} \times \eta'}$ with respect to the measure $x^{- (n + 1)} d V_\Leb$~.~$\Box$ \\

{\it Now we prove theorem \ref{main} . \\

Let $\rho \in \{0, \dots, r\}$ and $k \geq 2 n - \rho$ . Since $\vol \Gamma \backslash G < \infty$ it suffices to show that $f \in sM_k^{(\rho)}(\Gamma)$ and $\widetilde f \in L^1\rund{\Gamma \backslash G} \otimes \O\rund{\cz^{0| r}}$ imply \\
$\widetilde f \in L^\infty\rund{\Gamma \backslash G} \otimes \O\rund{\cz^{0| r}}$~. So let $f \in sM_k^{(\rho)}(\Gamma)$ such that \\
$\widetilde f \in L^1\rund{\Gamma \backslash G} \otimes \O\rund{\cz^{0| r}}$~. Let $g \in \Xi$ . By corollary \ref{L^s} it is even enough to show that
$l_g \rund{\widetilde f} \in L^\infty\rund{\eta A_{> t_0} K} \otimes \O\rund{\cz^{0| r}}$ , where $l_g \rund{\widetilde f}$ again denotes the left translation of $\widetilde f$ by the group element $g \in G$ . Let

\[
\left.f|_g \right|_{R^{- 1}} = \sum_{I \in \wp(r) \, , \, \abs{I} = \rho} q_I \bthet^I \, ,
\]

all $q_I \in \O(H)$ , $I \in \wp(r)$ , $\abs{I} = \rho$ . Then

\[
f|_g = \sum_{I \in \wp(r) \, , \, \abs{I} = \rho} q_I\rund{R \bw} \bzeta^I j\rund{R, \bw}^{k + \rho} \, .
\]

Since by corollary \ref{L^s} $\widetilde f \in L^1 \rund{\eta A_{> t_0} K} \otimes \O\rund{\cz^{0| r}}$ we conclude that

\[
q_I\rund{R \bz} j\rund{R, \bz}^{k + \rho} \Delta\rund{\bz, \bz}^{\frac{k + \rho}{2}} \in L^1\rund{\eta A_{> t_0} \b0}
\]

with respect to the $G$-invariant measure on $B$ or equivalently $q_I \Delta'\rund{\bw, \bw}^{\frac{k + \rho}{2}} \in L^1\rund{R \eta A_{> t_0} \b0}$ for all
$I \in \wp(r)$ , $\abs{I} = \rho$ , with respect to the $R G R^{- 1}$ -invariant measure on $H$ . So by theorem \ref{Satake Fourier} (iii) we see that
$q_I \Delta'\rund{\bw, \bw}^{\frac{k + \rho}{2}} \in L^\infty\rund{R \eta A_{> t_0} \b0}$ as well, or equivalently
$q_I\rund{R \bz} j\rund{R, \bz}^{k + \rho} \Delta\rund{\bz, \bz}^{\frac{k + \rho}{2}} \in L^\infty\rund{\eta A_{> t_0} \b0}$ for all $I \in \wp(r)$ , $\abs{I} = \rho$ . Therefore

\[
l_g\rund{\widetilde f} \in L^{\infty}\rund{\eta A_{> t_0} K} \otimes \O\rund{\cz^{0| r}} \, . \, \Box
\]

\end{document}